\documentclass[11pt]{article}
\usepackage[a4paper,top=3cm,bottom=2cm,left=2cm,right=2cm,marginparwidth=1.75cm]{geometry}

\usepackage[utf8]{inputenc} 
\usepackage[T1]{fontenc}
\usepackage[numbers]{natbib}
\usepackage[colorlinks=True,citecolor=blue,urlcolor=blue,pagebackref=true,backref=true]
{hyperref}
\renewcommand*{\backrefalt}[4]{%
    \ifcase #1 \footnotesize{(Not cited.)}%
    \or        \footnotesize{(Cited on page~#2.)}%
    \else      \footnotesize{(Cited on pages~#2.)}%
    \fi}

\usepackage[none]{hyphenat}
\usepackage{breakcites}
\usepackage[utf8]{inputenc} 
\usepackage[T1]{fontenc}    
\usepackage{hyperref}       
\hypersetup{
  colorlinks = true,
  urlcolor = blue,
  linkcolor = blue,
  citecolor = blue
}
\usepackage{booktabs}       
\usepackage{amsfonts}       
\usepackage{nicefrac}       
\usepackage{microtype}      
\usepackage{xcolor}         
\usepackage{caption}
\usepackage{subcaption}
\usepackage{float}
\usepackage{amsmath,amsthm}
\usepackage{amssymb}
\usepackage{mathtools}
\usepackage{natbib}
\usepackage{soul}
\usepackage{bm}
\usepackage{enumitem}
\usepackage{multirow}
\usepackage{wrapfig}
\usepackage{scalerel}
\allowdisplaybreaks
\usepackage{dsfont}
\usepackage[noend]{algpseudocode}
\usepackage{xcolor}
\usepackage{thmtools,thm-restate}
\usepackage{cleveref}
\usepackage{algorithm, algpseudocode}
\usepackage{url}

\title{Dual Control of Linear Systems from Bilinear Observations \\ with Belief Space Model Predictive Control}

\author{
\begin{tabular}{ccc}
Daniel Cao$^\ast$ & Beixi Du$^\ast$ & Andrew Lowitt$^\ast$ \\
Cornell & Cornell & Cornell \\
{\tt dyc33@cornell.edu} & {\tt bd469@cornell.edu} & {\tt ael254@cornell.edu} \\[2ex]
Sunmook Choi$^\ast$ & Sarah Dean & Yahya Sattar \\
Cornell & Cornell & Cornell \\
{\tt sc3377@cornell.edu} & {\tt sdean@cornell.edu} & {\tt ysattar@cornell.edu}
\end{tabular}
}

\date{}
\usepackage{enumitem}
\usepackage{caption,subcaption,multirow}
\usepackage{bm}

\renewcommand{\cite}{\citep}

\newcommand{\Ical}{\mathcal{I}}

\newcommand{\Ncal}{\mathcal{N}}

\newcommand{\E}{\operatorname{\mathbb{E}}}

\newcommand{\R}{\mathbb{R}}				

\newcommand{\dm}[2]
{
	\IfStrEq{#2}{1}{\R^{#1}}{\R^{#1 \x #2}}
}

\newcommand{\tr}{\textup{\textbf{tr}}} 			
\newcommand{\distas}{\overset{\text{i.i.d.}}{\sim}}

\newcommand{\T}{\top}
\newcommand*{\x}{\mathsf{x}\mskip1mu} 	

\newcommand{\splitatcommas}[1]{%
	\begingroup
	\begingroup\lccode`~=`, \lowercase{\endgroup
		\edef~{\mathchar\the\mathcode`, \penalty0 \noexpand\hspace{0pt plus .1em}}%
	}\mathcode`,="8000 #1%
	\endgroup
}

\usepackage{bbm}
\newcommand{\Item}[1]{%
	\ifx\relax#1\relax  \item \else \item[#1] \fi
	\abovedisplayskip=0pt\abovedisplayshortskip=0pt~\vspace*{-\baselineskip}
} 

\usepackage{booktabs}
\usepackage{multirow}

\usepackage{mathtools}
\makeatletter
\newcommand{\raisemath}[1]{\mathpalette{\raisem@th{#1}}}
\newcommand{\raisem@th}[3]{\raisebox{#1}{$#2#3$}}
\makeatother

\newcommand{\newcustomtheorem}[2]{%
	\newenvironment{#1}[1]
	{%
		\renewcommand\customgenericname{#2}%
		\renewcommand\theinnercustomgeneric{##1}%
		\innercustomgeneric
	}
	{\endinnercustomgeneric}
}
\newcustomtheorem{customtheorem}{Theorem}
\newcustomtheorem{customlemma}{Lemma}
\newcustomtheorem{customproposition}{Proposition}

\newcommand{\vertiii}[1]{{\vert\kern-0.25ex\vert\kern-0.25ex\vert #1 \vert\kern-0.25ex\vert\kern-0.25ex\vert}}

\newcommand{\beq}{\begin{equation}}
	
	\newcommand{\eeq}{\end{equation}}

\definecolor{emmanuel}{RGB}{255,127,0}

\numberwithin{equation}{section} 

\usepackage[]{mdframed}

\begin{document}

\maketitle

\begingroup
\renewcommand{\thefootnote}{\fnsymbol{footnote}}
\footnotetext[1]{Equal contribution.}
\endgroup



\begin{abstract}
We study finite-horizon quadratic control of linear systems with bilinear observations, in which the control input affects not only the state dynamics but also the partial observations of the state. In this setting, the separation principle can fail because control inputs influence the future quality of state estimates. State estimation requires an input-dependent Kalman filter whose gain and error covariance evolve as functions of the control inputs. To address this challenge, we propose a belief-space model predictive control (\texttt{B-MPC}) method that plans directly over both the estimated state and its error covariance.
In particular, \texttt{B-MPC} plans with a deterministic surrogate of the belief evolution defined by the input-dependent Kalman filter. Through numerical experiments in two synthetic settings, we show that \texttt{B-MPC} can outperform both the separation-principle controller and its MPC variant in favorable regimes, and that these gains are accompanied by lower estimation covariance and more uncertainty-aware action choices.
\end{abstract}

\section{Introduction}

Linear dynamical systems with bilinear observations (LDS-BO) are a simple example of the \emph{observer effect}, in which measurements about a system's state also alter that state.
Observer effects are present in applications ranging from quantum systems \cite{pardalos2010optimization} to electronic circuits \cite{liu2024system}. 
From the measurement perspective, estimating an unknown state is complicated by the fact that the measurement mechanism may itself influence the system.
Motivated by this measurement question, several recent papers study system identification of unknown dynamics in the LDS-BO setting, either through estimating the input-output behavior \cite{choi2025explore}, a state space realization \cite{liu2024system}, or both \cite{sattar2025learning}.

The presence of observer effects also complicates control design,
even when the dynamics are known.
This is because control inputs play a \emph{dual role} of both affecting and observing the state.
This trade-off is evident in high-performance robotics applications such as drone racing, where accounting for uncertainty during planning can be critical for performance \cite{kaufmann2023champion}.
The dual role of control inputs can be made explicit in the LDS-BO setting with Gaussian noise:
the state is linearly affected through the control matrix, while the estimation covariance is nonlinearly affected through the (input-dependent) measurement matrix.

Perhaps the simplest optimal control problem that captures the dual control phenomena is quadratic control of LDS-BO with Gaussian noise (BO-LQG). 
This is a slight variation on the classical linear quadratic Gaussian (LQG) problem,
which is solved by the separation principle.
Eschewing uncertainty, the optimal LQR controller is linear state feedback treating the state estimate as if it were true.
In contrast, separation-principle controller is not optimal for BO-LQG.
A recent work \cite{sattar2025sub} shows that the optimal controller is not linear in the estimated state and that the separation principle controller may locally maximize the cost.

Solving the BO-LQG problem is nontrivial.
It is a specific instance of a partially observed Markov decision process,
which are known to be computationally hard in general \cite{papadimitriou1987complexity},
largely due to the complexity of reasoning about long sequences of actions and observations \cite{kaelbling1998planning}.
However, the BO-LDS setting admits a finite dimensional \emph{belief state}: the distribution of the unobserved state is completely characterized by the state estimation mean and covariance.
Belief space planning is a common strategy in robotics to handle practical issues like path planning under perception uncertainties \cite{nagami2024state,zheng2022belief,zinage2023optimal}.
We take inspiration from this perspective, which allows us to reformulate the original optimal control problem into a belief space control problem.
This transforms the challenge from reasoning about partial observation to synthesizing a controller for nonlinear dynamics, which arise from the nonlinear covariance evolution.

A popular approach for addressing nonlinear optimal control problems is model predictive control (MPC) \cite{grune2016nonlinear, mayne2014model,rawlings2024model}. 
Instead of synthesizing a controller directly, an open loop planning problem is repeatedly solved, and each time only the first input is applied.
This strategy is usually an approximation: the planning horizon is shortened and surrogate deterministic dynamics may be used (though a variety of approaches exist~\cite{rawlings2024model}).
Indeed, MPC strategies have also been used for control of bilinear dynamical systems \cite{kafash2022model,xie2025bilinear}, however, these methods differ from our setting in that they consider the full state observation with
the bilinearity appearing in the state dynamics rather than the observation equation.

\emph{Contributions:} In this paper, we make the following contributions. First, we formulate linear quadratic control of linear systems with bilinear observations as a belief-space control problem, where the belief state is given by the Kalman filter state estimate and error covariance. This formulation makes the dual role of control explicit in the planning problem. Second, we propose a belief-space model predictive control method, \texttt{B-MPC}, that plans over a deterministic surrogate of the belief dynamics induced by the input-dependent Kalman filter.
Lastly, through numerical experiments on two synthetic systems, we show that \texttt{B-MPC} can substantially outperform the separation-principle controller and its MPC variant in favorable parameter regimes.
The experiments further show that these gains are associated with lower estimation covariance and improved observability, supporting the dual-control interpretation of the problem.

\emph{Organization:} The rest of the paper is organized as follows. In \S\ref{sec:problem-setting}, we formally introduce BO-LQG, and discuss the sub-optimality of the separation-principle controller \cite{sattar2025sub}. In \S\ref{sec:mpc-formulation}, we develop the belief-space MPC formulation and present the \texttt{B-MPC} algorithm.
In \S\ref{sec:experiments-setting-performance}, we describe the experimental setup and evaluate the performance of \texttt{B-MPC} against baseline controllers.
In \S\ref{sec:experiments-interpretation-discussion}, we further interpret the empirical behavior of \texttt{B-MPC} relative to the baselines through analyses of state estimation and control behavior. Lastly, \S\ref{sec:conclusion} concludes the paper.

\emph{Notation:} For a symmetric matrix $M$, we write $M \succeq 0$ (or $M \succ 0$) for positive semidefinite (or positive definite). For an arbitrary square matrix $N$, we write $\tr(N)$ for the trace and $\rho(N)$ for the spectral radius. We use $I_n$ to denote the $n\times n$ identity matrix. For a vector $u$, $(u)_i$ denotes the $i$-th element. We use the shorthand for sequences $u_{t:t+k}=(u_t,...,u_{t+k})$.

\section{Problem Settings} \label{sec:problem-setting}

In this paper, we consider the following linear dynamical system with bilinear observations:
\begin{equation}
\begin{aligned} \label{eqn: dynamics update}
    	x_{t+1} &= A x_t + Bu_t + w_{t}, \\
    y_t &= \big(C_0 + \sum_{k=1}^p(u_t)_kC_k\big) x_t + z_t.
\end{aligned}
\end{equation}
Here, $x_t {\in}\R^n$, $u_t {\in} \R^p$, and $y_t {\in} \R^m$ are the state, input, and the output at time $t$ with the initial state distribution $x_0 {\sim} \Ncal(\hat x_0 , \Sigma_0)$.
The matrices $A \in \R^{n\times n}$ and $B \in \R^{n\times p}$ define the linear dynamics, and $C_0,C_1,\dots,C_p \in \R^{m \times n}$ define the partial observation model. 
In this work, we assume that the matrices $A,B,C_0,C_1,\dots,C_p$ are known.
We assume that the process noise $w_t \in \R^n$ and the measurement noise $z_t \in \R^m$ are independent Gaussian processes, i.e., $w_t {\distas} \Ncal(0,\Sigma_w)$ and $z_t {\distas} \Ncal(0,\Sigma_z)$.
Given an input vector $u_t {\in} \R^p$, the matrix $C(u_t)$ denotes the input dependent observation matrix, i.e., $C(u_t) {=} C_0 {+} \sum_{k=1}^p (u_t)_k C_k$. Letting $\Theta = (A,B,C_0,C_1,\dots,C_p)$, we use LDS-BO($\Theta$) to refer to the state-space representation~\eqref{eqn: dynamics update}. 
Under the partially observed system, we want to solve the following finite-horizon optimal control problem, which we refer to as BO-LQG:
\begin{equation}
\begin{aligned}\label{eqn:bilinear LQG}
    \min_{u_{0:T-1}} &\E \left[ x_T^\T Q_T x_T + \sum_{t=0}^{T-1} \left( x_t^\top Q x_t + u_t^\top R u_t \right) \right] \quad \text{s.t.} \quad \text{LDS-BO}(\Theta),
\end{aligned}
\end{equation}
where $Q_T,Q \succeq 0$ are known state-cost matrices and $R \succ 0$ is a known input-cost matrix.

One can observe that, if the matrices $C_1=\cdots=C_p=0$, it reduces to the classical LQG setting. 
In classical LQG, the separation principle states that state estimation and control design can be solved independently: the state is estimated via Kalman filtering \cite{kalman1960new}, and the control input is obtained by applying the optimal linear feedback gain to the state estimate.

A natural baseline in BO-LQG is the separation-principle controller, i.e.,
$u_t^\texttt{Sep} = L_t \hat{x}_{t|t-1}$
where $L_t$ is the finite-horizon state-feedback gain computed from the standard Riccati recursion.
Here, $\hat{x}_{t|t-1}$ denotes the state estimate, i.e., the conditional mean of $x_t$ given the information $\Ical_t = \{u_0,\dots,u_{t-1},y_0,\dots,y_{t-1}\}$ available at time $t$.
Since the observation model depends on the input through $C(u_t)$, this estimate is obtained via an input-dependent Kalman filter.
In particular, the state estimate and the corresponding prediction error covariance evolve according to
\begin{equation} 
\begin{aligned} \label{eqn:kalman-recursion}
    \hat x_{t+1|t} &= A \hat x_{t|t-1} + B u_t - L(u_t)\left(y_t - C(u_t) \hat x_{t|t-1} \right) \\
    \Sigma_{t+1|t} &= A \Sigma_{t|t-1} A^\top + L(u_t)C(u_t) \Sigma_{t|t- 1} A^\top + \Sigma_w
\end{aligned}
\end{equation}
where $L(u_t) = -A\Sigma_{t|t{-}1}C(u_t)^{\top}
(C(u_t) \Sigma_{t|t{-}1}C(u_t)^{\top} + \Sigma_z)^{-1}$ is the input-dependent Kalman gain. 

Prior work \cite{sattar2025sub} showed that the separation principle does not hold in BO-LQG. 
In particular, it showed that the optimal controller for \eqref{eqn:bilinear LQG} is not affine in the estimated state, and provided a setting in which the separation principle controller locally maximizes the cost function.

The input-dependent Kalman filter also helps explain this suboptimality.
Since the observation matrix $C(u_t)$ depends on the applied input, the quality of the state estimate depends on the control sequence itself. 
If $C(u_t)$ is weakly informative, then the resulting observation provides little information about the latent state, and the error covariance $\Sigma_{t|t-1}$ remains large. 
Thus, an input that appears favorable for immediate regulation may simultaneously degrade future state estimation. 
The separation-principle controller does not explicitly account for this effect of control on future uncertainty, which is why it can be suboptimal in the bilinear observation setting. 
The key distinction from classical LQG is that the same input $u_t$ appears both in the state transition and in the observation matrix. Thus, choosing $u_t$ affects not only the next state distribution but also the information obtained about the current state.


\section{Model Predictive Control Formulation}
\label{sec:mpc-formulation}

We now formulate model predictive control (MPC) for the
bilinear observation system LDS-BO$(\Theta)$ and
a finite planning horizon $H \ge 1$.
A natural formulation of the MPC problem is to minimize the following objective over $u_{t:t+H-1}$: 
\begin{align} \label{eqn:MPC-information}
    \E \! \Bigg[ x_{t+H}^\top Q_T x_{t+H} {+} \sum_{\tau=t}^{t+H-1} \! \left(x_\tau^\top Q x_\tau {+} u_\tau^\top R u_\tau \right) \,\Big|\, \Ical_t
\Bigg]
\end{align}
Let $J_t^{\mathrm{MPC},H}$ be the minimum of \eqref{eqn:MPC-information}.
In this section, we formulate an MPC problem in the belief space, which is equivalent to minimizing \eqref{eqn:MPC-information}, induced by the input-dependent Kalman filtering~\eqref{eqn:kalman-recursion}.

\subsection{Belief-State Cost and Finite-Horizon MPC Objective}

We define the belief state as the mean and covariance of the posterior distribution of the state given past observations and controls:
\begin{align*}
    b_t := (\hat x_{t|t-1}, \Sigma_{t|t-1})
\end{align*}
where $\hat x_{t|t-1}$ and $\Sigma_{t|t-1}$ are given by the input-dependent Kalman filter in \eqref{eqn:kalman-recursion}.
Note that the Kalman filter provides the exact Bayesian posterior due to the Gaussian noise and linear dynamics.
The posterior distribution of the state $x_t$ given $\Ical_t$ remains Gaussian, and is fully characterized by the mean and covariance in $b_t$.

Then, the original finite-horizon quadratic cost \eqref{eqn:bilinear LQG} can be
rewritten in terms of the belief state using
\begin{align*}
    \E \big[ x_t^\top Q x_t \mid \Ical_t \big] &= \hat x_{t|t-1}^\top Q \hat x_{t|t-1} + \tr(Q\Sigma_{t|t-1}) \\
    \E \big[x_T^\top Q_T x_T \mid \Ical_T \big] &= \hat x_{T|T-1}^\top Q_T \hat x_{T|T-1} + \tr(Q_T\Sigma_{T|T-1}).
\end{align*}
This motivates the following belief-space stage and terminal costs:
\begin{align}
\ell(b_t,u_t)
&:= \hat x_{t|t-1}^\top Q \hat x_{t|t-1}
   + \tr(Q \Sigma_{t|t-1})
   + u_t^\top R u_t, \notag \\
\phi(b_T)
&:= \hat x_{T|T-1}^\top Q_T \hat x_{T|T-1}
   + \tr(Q_T \Sigma_{T|T-1}). \label{eq:mpc-stage-terminal-cost} 
\end{align}

From these costs, we construct the finite-horizon MPC objective with planning horizon $H$ in the belief space which is equivalent to \eqref{eqn:MPC-information}.
Given the current belief $b_t$, we have that $J_t^{\mathrm{MPC},H} $ is equal to
\begin{align} \label{eq:mpc-belief-objective}
\min_{u_{t:t{+}H{-}1}}\E \! \Bigg[ \phi(b_{t+H}) +
\sum_{\tau=t}^{t+H-1} \ell(b_\tau,u_\tau)
\;\Big|\; b_t
\Bigg]
\end{align}
where the belief trajectory $(b_\tau)_{\tau=t}^{t+H}$ is generated by the
input-dependent Kalman filter \eqref{eqn:kalman-recursion} under the chosen open-loop control sequence
$(u_\tau)_{\tau=t}^{t+H-1}$ and the bilinear observation model
$y_\tau {=} C(u_\tau)x_\tau {+} z_\tau$.
When $H {=} T{-}t$, $J_t^{\mathrm{MPC},H}$ coincides with the optimal cost-to-go of the
full dynamic program. 

\subsection{Deterministic MPC Planning Objective}

In practice, planning with the objective \eqref{eq:mpc-belief-objective} is infeasible because the belief state update is stochastic and depends on the bilinear observation $y_\tau {=} C(u_\tau) x_\tau {+} z_\tau$.
Instead, we consider its deterministic surrogate based on the available information: $y_\tau {=} C(u_\tau)\hat x_\tau$.
Therefore, given an input $u_t$, we consider the following deterministic update for the planning:
\begin{equation} \label{eqn:det-belief-update}
    \begin{aligned}
        \bar x_{t+1} &= A \bar x_t + Bu_t \\
        \bar \Sigma_{t+1} &= A \bar\Sigma_t A^\top +\bar L(u_t)C(u_t) \bar\Sigma_{t}A^\top + \Sigma_w
    \end{aligned}
\end{equation}
where 
$\bar L(u_t) = -A\bar\Sigma_tC(u_t)^\top\big( C(u_t) \bar\Sigma_{t} C(u_t)^\top + \Sigma_z\big)^{-1}$. 
Under the nominal observation $y_\tau {=} C(u_\tau)\bar x_\tau$ the innovation term is zero, so the deterministic mean update reduces to the open-loop prediction $\bar x_{\tau+1} = A \bar x_\tau + Bu_\tau$ while the covariance still evolves according to the input-dependent Kalman covariance recursion.
Let $\bar b_t = (\bar x_t, \bar \Sigma_t)$, then \eqref{eqn:det-belief-update} defines a function $F$ such that
\begin{align} \label{eqn:belief-update-F}
    \bar b_{t+1} = F(\bar b_t, u_t).
\end{align}
Thus, the MPC objective corresponding to the deterministic planning is the following:
\begin{align} \label{eqn:B-MPC-objective}
    \bar J_t^{\mathrm{MPC},H} = \min_{u_{t:t+H-1}} \phi(\bar b_{t+H}) + \sum_{\tau = t}^{t+H-1} \ell(\bar b_\tau, u_\tau) 
\end{align}
where the trajectory $(\bar b_\tau)_{\tau=t}^{t+H}$ is generated by \eqref{eqn:belief-update-F} starting from $\bar b_t = (\hat x_{t|t-1},\Sigma_{t|t-1})$, the true belief state.
We emphasize that \eqref{eqn:B-MPC-objective} is an approximation to the stochastic belief-space objective \eqref{eq:mpc-belief-objective}.
Replacing future observations by their nominal values is a standard certainty-equivalent approach in nonlinear MPC. 
Our contribution is a tractable receding-horizon approximation for BO-LQG rather than an exact solution or a closed-loop stability guarantee, though it does account for the dual role of control inputs.
In Section~\ref{sec:experiments-interpretation-discussion} and \ref{sec:experiments-setting-performance}, we demonstrate its performance against separation-principle baselines and provide experiments that support the dual-control interpretation of the belief-space controller.
Unlike Sep-MPC, the B-MPC objective contains the trace terms  $\text{tr}(Q\Sigma_{t|t-1})$ and $\text{tr}(Q_T \Sigma_{t|t-1})$, and the covariance trajectory depends on the planned inputs through $C(u)$. Therefore, B-MPC can select inputs that are suboptimal for immediate state regulation but beneficial for future estimation. This is the mechanism by which the controller captures the dual role of control.

\subsection{Belief-Space MPC Algorithm} \label{sec:b-mpc-algorithm}

We state the receding-horizon MPC algorithm in belief space for any total horizon $T$ and planning horizon $H \ge 1$.

Given the system LDS-BO$(\Theta)$, cost
matrices $(Q,Q_T,R)$, noise covariances $(\Sigma_w,\Sigma_z)$, a look-ahead
window $H \ge 1$, we define a belief-space
MPC controller as follows. 
The controller maintains a belief state $b_t = (\hat x_{t|t-1},\Sigma_{t|t-1})$ at each time $t$, where $\hat x_{t|t-1}$ and $\Sigma_{t|t-1}$
are the mean and covariance of $x_t$ conditioned on the information $\Ical_t$. Starting from
the prior belief $b_0 = (\hat x_{0},\Sigma_{0})$, the mean and covariance of the initial state distribution, the controller repeats the
following steps for $t=0,1,\dots,T-1$:
\begin{enumerate}
    \item \textbf{Planning step.} Given the current belief $b_t$, the controller solves the
    finite-horizon optimal control problem
    \begin{equation*}
        \min_{u_{t:t+H-1}}
        \phi(\bar b_{t+H}) + \sum_{\tau = t}^{t+H-1} \ell(\bar b_\tau, u_\tau)
    \end{equation*}
    where the stage cost $\ell$ and terminal cost $\phi_{t+H}$ are defined in
    \eqref{eq:mpc-stage-terminal-cost}, and the belief trajectory
    $(\bar b_\tau)_{\tau=t}^{t+H}$ is generated by \eqref{eqn:belief-update-F}.
    
    \item \textbf{Control.} Let $u_{t:t+H-1}^\star$ denote
    an optimal solution of the above problem. The MPC controller applies only the first
    input,
    \[
        u_t := u_t^\star,
    \]
    to the true system.
    \item \textbf{Belief update.} After applying $u_t$, the controller observes
    $y_t$ and updates the belief to
    \[
        b_{t+1} = (\hat x_{t+1|t},\Sigma_{t+1|t}),
    \]
    using the bilinear Kalman filtering recursion~\eqref{eqn:kalman-recursion}.
\end{enumerate}
We remark that due to the belief state covariance update, the planning problem is nonconvex. However, it is differentiable. This motivates the use of gradient-based local optimization rather than dynamic programming or Riccati-style recursions.
We present further discussion on numerical optimization in the experiments section.
The algorithm for the belief-space MPC controller, called \texttt{B-MPC}, is summarized in Algorithm~\ref{alg:B-MPC}.
    
\begin{algorithm}[ht]
\caption{\texttt{B-MPC}} \label{alg:B-MPC}
\begin{algorithmic}[1]
\Require Trajectory length $T$, planning horizon $H$, mean $\hat x_0$ and covariance $\Sigma_0$ of the initial state distribution

\State Initialize $\hat x_{0|-1} = \hat x_0$ and $\Sigma_{0|-1} = \Sigma_0$.
\State Construct the belief state $b_0 = (\hat x_{0|-1}, \Sigma_{0|-1})$.
\For{$t=0,1,2,\dots,T-1$}
    \State From $b_t$, define a function $J$ by
    \begin{align*}
        J(u_t,\dots,u_{t{+}H{-}1}, b_t) = \sum_{\tau = t}^{t+H-1} \ell(\bar b_\tau, u_\tau) + \phi(\bar b_{t+H})
    \end{align*} 
    where $(\bar b_\tau)_{\tau=t}^{t+H}$ is generated by \eqref{eqn:belief-update-F} given the inputs $u_t,\dots,u_{t+H-1}$, starting from $\bar b_t = b_t$.
    \State Minimize $J$ over $(u_t,\dots,u_{t+H-1})$ via L-BFGS optimization, and let  $(u_t^\star,\dots,u_{t+H-1}^\star)$ denote a minimizer.
    \State Apply $u_t^\star$ in the true system, and observe $y_t$.
    \State Update the belief state $b_{t+1} = (\hat x_{t+1|t}, \Sigma_{t+1|t})$ via \eqref{eqn:kalman-recursion} using $(u_t^\star,y_t)$.
\EndFor
\end{algorithmic}
\end{algorithm}

\section{Experiments: Setting \& Performance} \label{sec:experiments-setting-performance}

In this section, we demonstrate the performance of our method by comparing it with two baseline controllers on two benchmark systems. 
We first define the benchmark systems and the controllers in Section~\ref{sec:bilinear-system-description} and Section~\ref{sec:controllers}.
In Section~\ref{sec:results-best}, we present settings in which \texttt{B-MPC}, our proposed belief-space MPC controller, achieves better performance than the baselines, and we analyze the resulting control behavior.

\subsection{Bilinear observation systems} \label{sec:bilinear-system-description}

We consider two bilinear observation systems: a multi-block double integrator system and a randomly generated system. In both systems, the initial state distribution is set to $\Ncal(0,I)$. We describe the system parameters in the following.

\subsubsection{Random bilinear observation system} 

We sample a random bilinear observation system with state dimension $n {=} 6$, input dimension $p {=} 3$, and output dimension $m {=} 3$. Each entry of $A$ is drawn independently from $\mathcal{N}(0,1)$, and the matrix $A$ is rescaled so that $\rho(A)$ matches the target spectral radius. We draw $B$ entrywise from i.i.d. $\mathcal{N}(0, 1/n)$, the observation matrix $C_0$ entrywise from $\mathcal{N}(0, c_0^2 / m)$, and each matrix $C_k$, $k=1,\dots,p$, entrywise from $\mathcal{N}(0, 1 / m)$. 
The noise covariances are chosen to $\Sigma_w = \sigma_w^2I_n$ and $\Sigma_z = \sigma_z^2I_m$, and the cost matrices are $Q = Q_T = I_n$ and $R = \mathrm{R}_{\mathrm{scale}} \cdot I_p$. 

\subsubsection{Multi-block double integrator style system} \label{sec:multi-block-double-integrator}

The multi-block double integrator system consists of three identical, dynamically decoupled blocks. 
The matrices $A$ and $B$ are defined to be
\begin{align*}
    A = \begin{bmatrix}
        \rho & h & 0 & 0 & 0 & 0 \\
        0 & \rho & 0 & 0 & 0 & 0 \\
        0 & 0 & \rho & h & 0 & 0 \\
        0 & 0 & 0 & \rho & 0 & 0 \\
        0 & 0 & 0 & 0 & \rho & h \\
        0 & 0 & 0 & 0 & 0 & \rho 
    \end{bmatrix}, \quad 
    B = \begin{bmatrix}
        0 & 0 & 0 \\ 
        h & 0 & 0 \\ 
        0 & 0 & 0 \\ 
        0 & h & 0 \\ 
        0 & 0 & 0 \\ 
        0 & 0 & h 
    \end{bmatrix}
\end{align*}
where $\rho$ is set to the target spectral radius and $h=0.3$ is the discretization step size.
The observation matrix \(C_0 \in \mathbb{R}^{3 \times 6}\) is defined by $(C_0)_{i,\,2i-1}=c_0$, $i=1,2,3$, with all other entries equal to zero. Equivalently,
\[
C_0=
\begin{bmatrix}
c_0 & 0 & 0 & 0 & 0 & 0\\
0 & 0 & c_0 & 0 & 0 & 0\\
0 & 0 & 0 & 0 & c_0 & 0
\end{bmatrix}.
\]
For $k{=}1,2,3$, the matrix $C_k {\in} \mathbb{R}^{3 \times 6}$ is defined so that its $(k,\,2k{-}1)$-entry is $c_1$ and all other entries are zero. Throughout the experiments, we set $c_1{=}3$.
As for the random system, 
we consider isotropic noise processes, i.e., $\Sigma_w = \sigma_w^2I_6$ and $\Sigma_z = \sigma_z^2I_3$. For the quadratic costs, we set $Q=Q_T=I_6$ and $R=\mathrm{R}_{\rm scale} \cdot I_3$.

\subsection{Controllers} \label{sec:controllers}

We consider three controllers to evaluate their performance on the benchmark systems.
The first is the sub-optimal separation principle controller, referred to as \texttt{Sep}, defined as follows:
\begin{equation}
    \begin{aligned} \label{eqn:sep-controller}
        u_t^{\texttt{Sep}} &= L_t \hat x_{t|t-1} \\
        \text{where} \quad L_t &= -(B^\top K_{t+1} B + R)^{-1}B^\top K_{t+1}A.
    \end{aligned}
\end{equation}
Here, the matrix $K_t$ is computed recursively through Riccati equation, starting from $K_T = Q_T$:
\begin{equation} 
    \begin{aligned} \label{eqn:sep-Riccati}
        K_t &= A^\top K_{t+1}A - P_t + Q, \quad \text{where} \\
        P_t &= A^\top K_{t+1}B(B^\top K_{t+1} B + R)^{-1}B^\top K_{t+1}A.
    \end{aligned}
\end{equation}
The second is the controller \texttt{Sep-MPC}, which is an MPC version of separation principle controller. At each time $t$, it solves the finite-horizon deterministic linear-quadratic problem:
\begin{align*}
    \min_{u_{0:H-1}} & \bar x_H^\top Q_T \bar x_H + \sum_{\tau=0}^{H-1} (\bar x_\tau^\top Q \bar x_\tau + u_\tau^\top Ru_\tau) \\
    \text{s.t. } & \bar x_{\tau+1} = A \bar x_\tau + B u_\tau, \quad \bar x_0 = \hat x_{t|t-1}.
\end{align*}
Equivalently, this problem is solved via the finite-horizon Riccati recursion
\begin{align*}
    K_t^{\rm MPC} &= A^\top K_{t+1}^{\rm MPC}A - P_t + Q, \quad \text{where} \\
    P_t &= A^\top K_{t+1}^{\rm MPC}B(B^\top K_{t+1}^{\rm MPC} B + R)^{-1}B^\top K_{t+1}^{\rm MPC}A
\end{align*}
that starts with $K_{t+H}^{\rm MPC} = Q_T$.
Then it applies the action $u_t = L_t^{\rm MPC} \hat{x}_{t|t-1}$ where 
\begin{align*}
    L_t^{\rm MPC} = -(B^\top K_{t+1}^{\rm MPC}B + R)^{-1} B^\top K_{t+1}^{\rm MPC}A.
\end{align*}

Finally, the third controller is
\texttt{B-MPC}, our proposed belief-space MPC controller. The algorithm is described in Section~\ref{sec:b-mpc-algorithm} and Algorithm~\ref{alg:B-MPC}.
At each planning step, \texttt{B-MPC} approximately solves the nonconvex finite-horizon optimization problem.
We use limited-memory BFGS (L-BFGS), a quasi-Newton optimization algorithm, implemented in PyTorch.
We initialize 
each entry of the $H$-length control sequence randomly i.i.d. from $\Ncal(0,1/H)$ at every time step, and the optimizer is run for 20 outer iterations with step size 0.8. We use a single initialization per planning step.


Unless noted otherwise, the experiments use $T {=} 300$ time steps, and the same sampled initial condition and noise realization across controllers within each trial.

\subsection{Improved performance with \texttt{B-MPC}}
\label{sec:results-best}

\begin{table}[ht]
    \centering
    \caption{System Parameter Configuration}
    \label{tab:sys-param-config}
    \begin{tabular}{lccccc}
        \hline
        \textbf{System} & $\rho(A)$ & $c_0$ & $\mathrm{R}_{\rm scale}$ & $\sigma_w$ & $\sigma_z$ \\
        \hline
        Random System & 0.95 & 0.01 & 1 & 0.1 & 0.1 \\
        Double Integrator & 0.95 & 0.01 & 1 & 0.1 & 1.0 \\
        \hline
    \end{tabular}
\end{table}

We first demonstrate that \texttt{B-MPC} can achieve substantial cost reductions over both \texttt{Sep} and \texttt{Sep-MPC}, for the parameter settings shown in Table~\ref{tab:sys-param-config}.
These parameters reflect a setting where the system is not very stable, the input-independent observation has small magnitude, and the state and input costs are comparable.

\begin{figure}[ht]
    \centering
    \begin{subfigure}{0.35\textwidth}
        \centering
        \includegraphics[width=\linewidth]{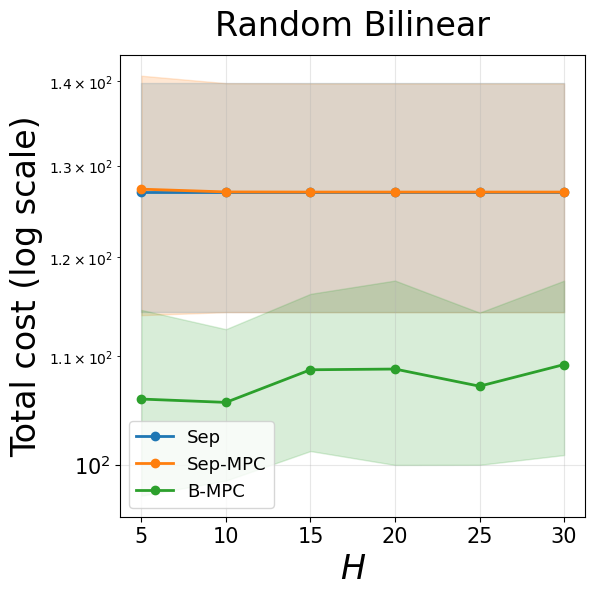}
        \label{fig:h-sweep-random}
    \end{subfigure}
    \hspace{10mm}
    \begin{subfigure}{0.35\textwidth}
        \centering
        \includegraphics[width=\linewidth]{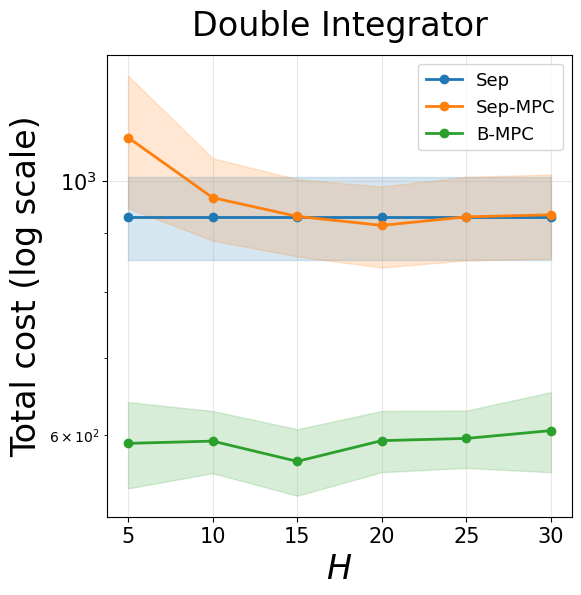}
        \label{fig:h-sweep-double}
    \end{subfigure}
    \caption{Total cost versus look-ahead horizon $H\in\{5,10,15,20,25,30\}$ for the three controllers, averaged over 10 trials. Shaded regions denote 95\% confidence intervals across trials.
    Left: random bilinear observation system. Right: multi-block double integrator system.}
    \label{fig:h-sweep}
\end{figure}

Under the parameter setting in Table~\ref{tab:sys-param-config}, Figure~\ref{fig:h-sweep} shows that \texttt{B-MPC} significantly outperforms the controllers \texttt{Sep} and \texttt{Sep-MPC} for all planning horizons $H \in \{5,10,15,20,25,30\}$.
For random bilinear system, the cost reduction from \texttt{B-MPC} is 16.9\% compared to the cost from \texttt{Sep} at $H{=}10$,
and for the double integrator system, the cost reduction is 38.8\% at $H{=}15$.
On the other hand, \texttt{Sep-MPC} performs worse than \texttt{Sep} for double integrator systems when the planning horizon is too low, and otherwise performs comparably.

To better understand the source of this improvement, we decompose the total cost into state and control components. The same parameters are used as in Table~\ref{tab:sys-param-config} with fixed look-ahead horizon $H$. We use $H=10$ and $H=15$ for random bilinear system and double integrator system, respectively.

\begin{table}[ht]
    \centering
    \caption{Cost decomposition for each controller in two systems. The average cost over 10 trials is shown.}
    \label{tab:cost-decomposition}
    \begin{tabular}{l@{\hspace{4pt}}l@{\hspace{5pt}}cccc}
        \hline
        \textbf{System} & $H$ & \textbf{Controller} & \textbf{State-cost} & \textbf{Input-cost} & \textbf{Total-cost} \\
        \hline
        \multirow{3}{*}{\shortstack{Random\\System}}
            & \multirow{3}{*}{10}
            & \texttt{Sep}     & 124.67 & 2.38 & 127.05 \\
            & 
            & \texttt{Sep-MPC} & 124.68 & 2.38 & 127.06 \\
            & 
            & \texttt{B-MPC}   & 87.06 & 18.57 & 105.64 \\
        \hline
        \multirow{3}{*}{\shortstack{Double\\Integrator}}
            & \multirow{3}{*}{15}
            & \texttt{Sep}     & 876.22 & 53.93 & 930.15 \\
            & 
            & \texttt{Sep-MPC} & 875.05 & 56.14 & 931.20 \\
            & 
            & \texttt{B-MPC}   & 422.64 & 146.30 & 568.94 \\
        \hline
    \end{tabular}
\end{table}

Table~\ref{tab:cost-decomposition} provides extra information in addition to Figure~\ref{fig:h-sweep}.
At a fixed look-ahead window $H$ for each system, it shows that \texttt{B-MPC} achieves substantially lower state cost compared to \texttt{Sep} and \texttt{Sep-MPC}, while it incurs slightly larger control cost.
This describes where \texttt{B-MPC} reduces the cost throughout the trajectory.


\section{Experiments: Interpretation \& Discussion}
\label{sec:experiments-interpretation-discussion}

In this section, we interpret the performance differences observed in Section~\ref{sec:results-best}. 
Our goal is to understand why \texttt{B-MPC} outperforms \texttt{Sep} and \texttt{Sep-MPC} in the favorable regime of Figure~\ref{fig:h-sweep}.
In particular, we examine the controllers from three complementary perspectives: estimation quality in Section~\ref{sec:kf-error}, the control actions selected along a trajectory in Section~\ref{sec:action-compare}, and how the discrepancy between \texttt{B-MPC} and \texttt{Sep-MPC} changes with the level of uncertainty in Section~\ref{sec:action-difference-synthetic}.
The results illustrate how the main advantage of \texttt{B-MPC} arises due to the fact that planning is carried out in the belief space.
For the experiments, we use the same parameter configuration as in Table~\ref{tab:sys-param-config} with fixed planning horizon $H{=}10$ and $H{=}15$ for random bilinear system and double integrator system, respectively.

\subsection{Estimation improves through covariance reduction} \label{sec:kf-error}

\begin{figure}[ht]
    \centering
    \includegraphics[width=0.75\linewidth]{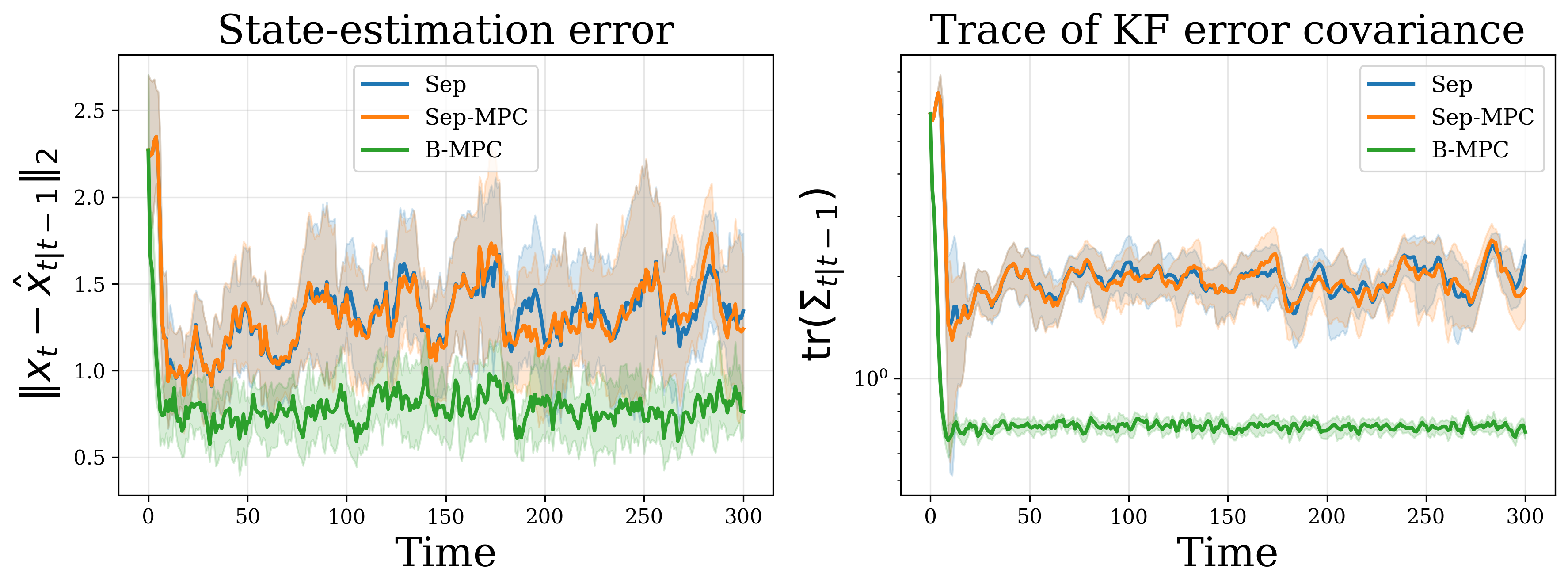}
    \caption{Kalman filter diagnostics on the multi-block double integrator system with $H{=}15$. Left: Kalman filter state estimation error $\|x_t - \hat{x}_{t|t-1}\|_2$. Right: the trace of Kalman filter estimation covariance $\tr(\Sigma_{t|t-1})$. Curves are averages over 10 trials under matched initial conditions and noise realizations. The shaded region denotes 95\% confidence intervals across trials.}
    \label{fig:kf-rollout}
\end{figure}

We examine trajectories of the multi-block double integrator style system, averaged over 10 trials.
Figure~\ref{fig:kf-rollout} shows that the performance gain of \texttt{B-MPC} is closely related to improved state estimation and lower error covariance.
On the left, the estimation error $\|x_t - \hat x_{t|t-1}\|_2$ is generally lower for \texttt{B-MPC} than for \texttt{Sep} and \texttt{Sep-MPC} over the trajectory. 
On the right, the panel shows an even stronger separation in the estimation covariance trace $\tr(\Sigma_{t|t-1})$, with \texttt{B-MPC} maintaining substantially smaller uncertainty throughout the trajectory.
Because the covariance update depends on the chosen input in the bilinear observation model, this gap indicates that \texttt{B-MPC} selects actions that make the system state easier to estimate.
Hence, Figure~\ref{fig:kf-rollout} provides an empirical illustration that \texttt{B-MPC} leverages the dual role of control, namely state regulation and improvement of future observability.

\subsection{Differences in control actions} \label{sec:action-compare}

We compare a \texttt{Sep-MPC} and \texttt{B-MPC} through a counterfactual action-matching experiment. First, we generate a single trajectory using \texttt{Sep-MPC} controller. This produces a sequence of belief states $(b_t)_{t=0}^{T-1}$ and inputs $(u_t^\texttt{Sep-MPC})_{t=0}^{T-1}$ along the trajectory. At each time step $t$, we solve \texttt{B-MPC} controller from the same belief state $b_t$, obtaining $u_t^\texttt{B-MPC}$.
By repeating this for all time steps, we obtain a sequence of counterfactual input sequence $(u_t^\texttt{B-MPC})_{t=0}^{T-1}$. 
The goal of this section is to compare the two input sequences.

\begin{figure}[!htbp]
    \centering
    \includegraphics[width=0.75\linewidth]{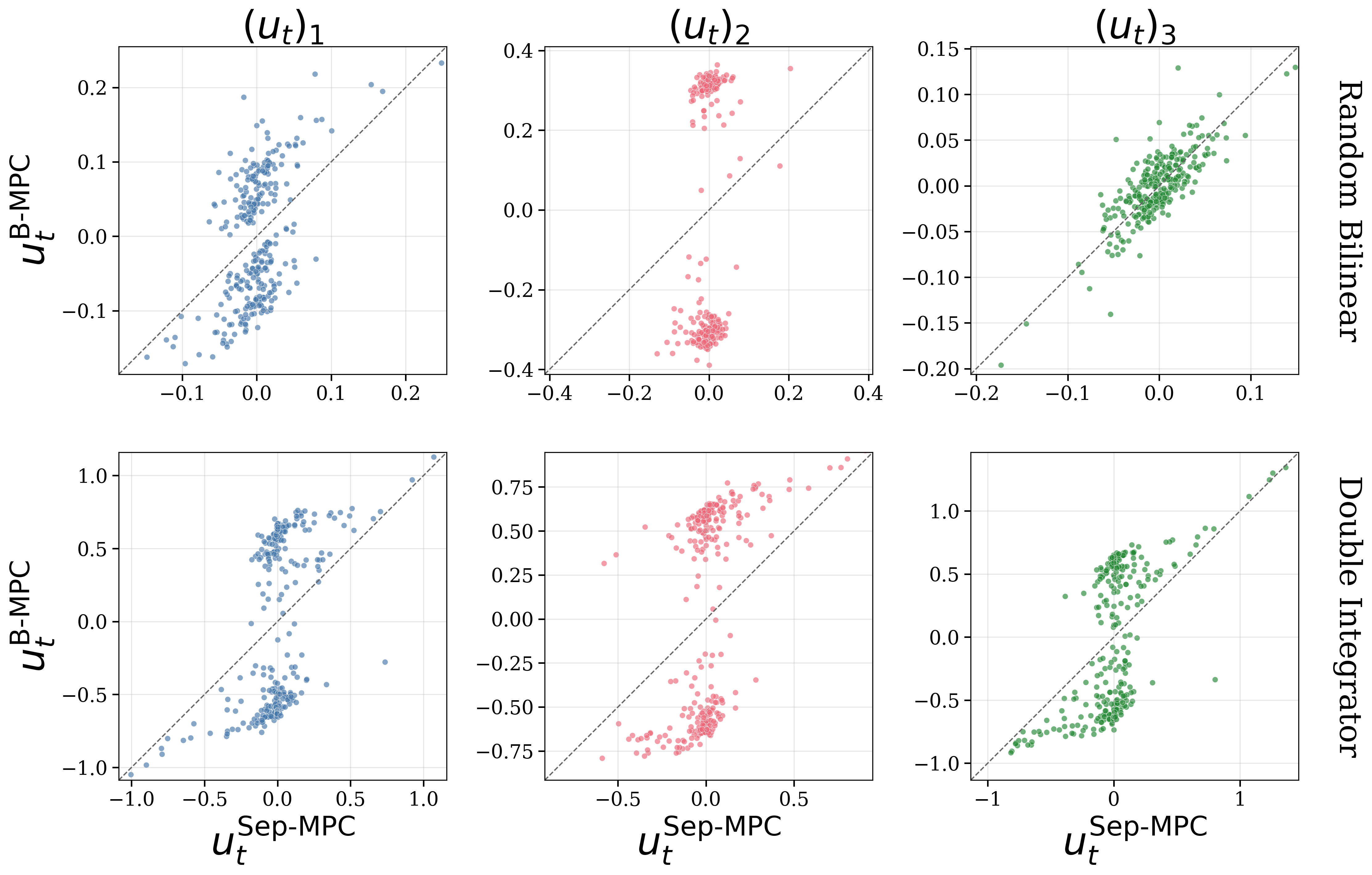}
    \caption{Counterfactual action comparison between \texttt{Sep-MPC} and \texttt{B-MPC}. A trajectory from \texttt{Sep-MPC} is generated first, and at each time step, \texttt{B-MPC} is solved from the same belief state. Each point shows the corresponding \texttt{Sep-MPC} and \texttt{B-MPC} action coordinates. The dashed line indicates equality, $u^{\texttt{B-MPC}} = u^{\texttt{Sep-MPC}}$.}
    \label{fig:input-scatter}
\end{figure}

Figure~\ref{fig:input-scatter} shows the result. 
From the left to right, each panel corresponds to the $i$-th element of the action.
From the two input sequences obtained above, 
we plot
\begin{align*}
    \big((u_t^{\texttt{Sep-MPC}})_i, (u_t^{\texttt{B-MPC}})_i \big) 
\end{align*}
so that the $x$-axis is the $i$-th coordinate of $u_t^{\texttt{Sep-MPC}}$ and $y$-axis is the corresponding coordinate of $u_t^\texttt{B-MPC}$. 
Points on the diagonal indicate agreement between the two controllers, while deviations from the diagonal indicate differences in the selected action coordinate.

In both systems, we can see that the actions from \texttt{Sep-MPC} are near zero while the ones from \texttt{B-MPC} are more spread. 
Interestingly, \texttt{B-MPC} tends to choose inputs that move away from zero, but not so far that the input penalty dominates.
This is consistent with dual-control interpretation of \texttt{B-MPC}: the controller explicitly trades off immediate regulation and future information.
Since $c_0{=}0.01$, near-zero actions make the input-dependent observation matrix $C(u_t)$ close to zero, which can substantially reduce observability. This results in worse state estimation and high uncertainty as in Figure~\ref{fig:kf-rollout}.

\subsection{Controller disagreement increases with uncertainty} \label{sec:action-difference-synthetic}

The scatter plots in Figure~\ref{fig:input-scatter} show that \texttt{B-MPC} and \texttt{Sep-MPC} disagree, but they do not reveal \emph{when} the disagreement is largest. Since \texttt{B-MPC} accounts for covariance evolution in its planning objective, we expect the action difference to grow with the level of state uncertainty. To test this directly, we construct synthetic beliefs: $10$ random state estimates $\hat{x} \sim \mathcal{N}(0, 0.25 I)$, each paired with $\Sigma = \alpha I$ for $20$ log-spaced values of $\alpha$ ($\tr(\Sigma) \in [0.06, 190]$). At each synthetic belief state, we solve both \texttt{B-MPC} and \texttt{Sep-MPC} and record $\|u^{\texttt{B-MPC}} - u^{\texttt{Sep-MPC}}\|_2$.

\begin{figure}[!htbp]
    \centering
    \includegraphics[width=0.75\linewidth]{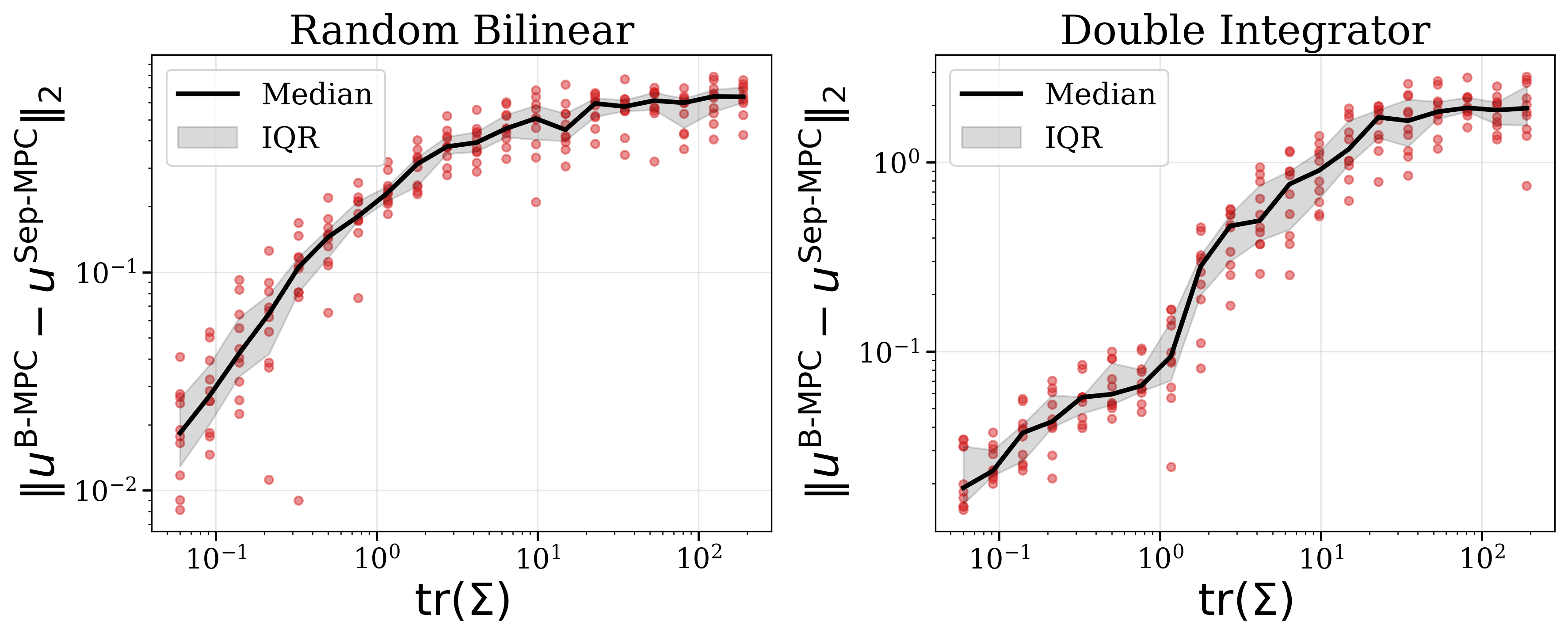}
    \caption{Action difference $\|u^{\texttt{B-MPC}} - u^{\texttt{Sep-MPC}}\|_2$ versus $\tr(\Sigma)$ for synthetic belief states. Results are computed from 10 sampled state estimates combined with 20 log-spaced covariance scales. The black line shows the median across synthetic beliefs, and the shaded region indicates the interquartile range.}
    \label{fig:action-gap-synthetic}
\end{figure}

Figure~\ref{fig:action-gap-synthetic} shows the result on a log-log scale. For both systems, the gap between \texttt{B-MPC} and \texttt{Sep-MPC} increases monotonically as $\tr(\Sigma)$ grows.
This confirms that \texttt{B-MPC} deviates from \texttt{Sep-MPC} when uncertainty is high. 
This behavior aligns with the fundamental trade-off of dual control, demonstrating that the belief-space controller explicitly weighs the long-term value of uncertainty reduction against the immediate penalty of cost minimization.

\subsection{Regimes where \texttt{B-MPC} provides limited benefit}
\label{sec:results-heatmap}

Separation principle remains a popular heuristic, and may perform reasonably in practice.
Indeed, we find that belief space planning does not always provide benefits:
\texttt{B-MPC} does not uniformly outperform the baselines. To map out the boundary, we compute the percentage cost reduction of \texttt{B-MPC} over \texttt{Sep} across a grid of $\mathrm{R}_{\mathrm{scale}} \in \{1, 10, 100\}$ and $c_0 \in \{0.01, 0.1, 1.0\}$, selecting for each combination the horizon $H \in \{5,10,15,20,25,30\}$ that minimizes \texttt{B-MPC} cost.

\begin{figure}[ht]
    \centering
    \includegraphics[width=0.75\linewidth]{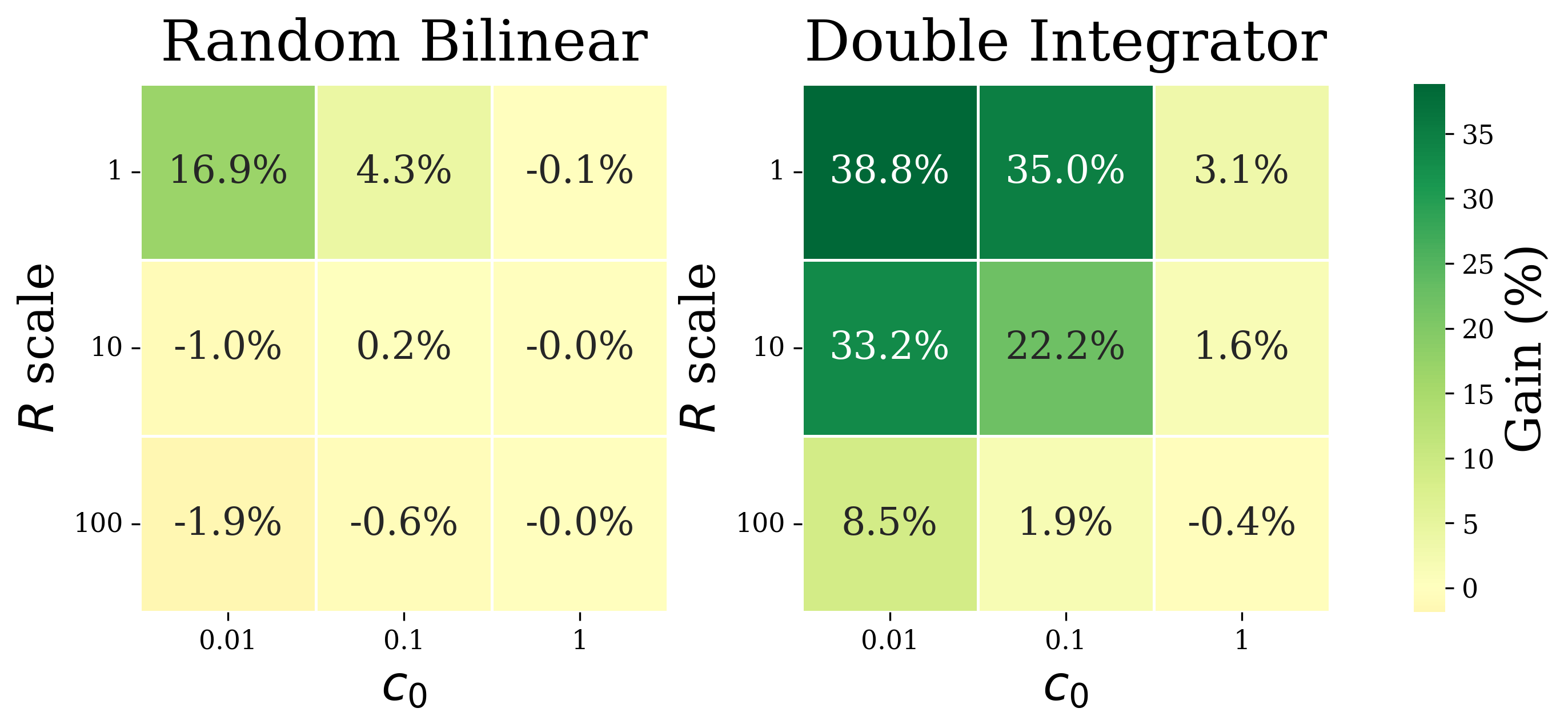}
    \caption{Percentage cost improvement of \texttt{B-MPC} over \texttt{Sep} across $\mathrm{R}_{\rm scale}$ and $c_0$ settings for the random bilinear and double-integrator systems, with $\rho(A)=0.95$. For each $(\mathrm{R}_{\rm scale}, c_0)$ pair, we select the horizon $H$ that minimizes the mean \texttt{B-MPC} cost over 10 trials, and report the relative gain  $100\times (J_\texttt{Sep} - J_\texttt{B-MPC})/J_\texttt{Sep}$ (\%). Each cell is annotated with the corresponding percentage. }
    \label{fig:heatmap-sep}
\end{figure}

In Figure~\ref{fig:heatmap-sep}, we can observe a trend that \texttt{B-MPC} tends to perform better as both $c_0$ and $\mathrm{R}_{\rm scale}$ are smaller.
If $c_0$ is small, then the input-dependent observation matrix $C(u_t)$ highly depends on the current input $u_t$.
Due to the quadratic input cost, \texttt{Sep} finds inputs that are close to zero, resulting in high chance of observability loss. 
On the other hand, as \texttt{B-MPC} accounts for state uncertainty through the belief space, it outperforms \texttt{Sep} when $\mathrm{R}_{\rm scale}{=}1$.
However, as $\mathrm{R}_{\rm scale}$ increases, the performance gap between \texttt{B-MPC} and \texttt{Sep} decreases, since large control penalties restrict the controller's ability to choose exploratory actions that reduce uncertainty.

\begin{table}[ht]
    \centering
    \caption{System Parameter Configuration}
    \label{tab:sys-param-config2}
    \begin{tabular}{lccccc}
        \hline
        \textbf{System} & $H$ & $c_0$ & $\mathrm{R}_{\rm scale}$ & $\sigma_w$ & $\sigma_z$ \\
        \hline
        Random System & 10 & 0.01 & 1 & 0.1 & 0.1 \\
        Double Integrator & 15 & 0.01 & 1 & 0.1 & 1.0 \\
        \hline
    \end{tabular}
\end{table}

For an additional view of the regime in which belief-space planning is less favorable, we vary the stability of the system. The system parameters are set according to Table~\ref{tab:sys-param-config2}.

\begin{figure}[ht]
    \centering
    \begin{subfigure}{0.35\textwidth}
        \centering
        \includegraphics[width=\linewidth]{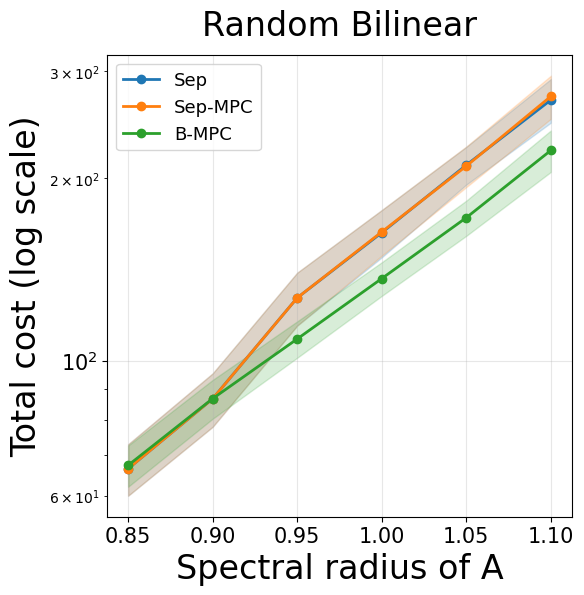}
        \label{fig:rho-sweep-random}
    \end{subfigure}
    \begin{subfigure}{0.35\textwidth}
        \centering
        \includegraphics[width=\linewidth]{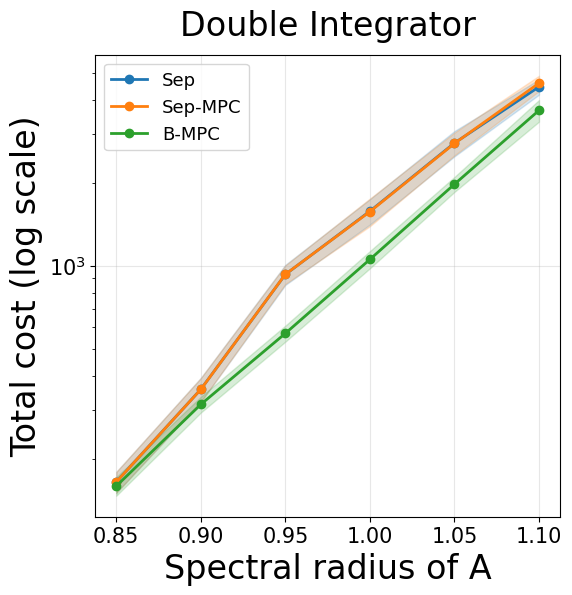}
        \label{fig:rho-sweep-double}
    \end{subfigure}
    \caption{Total cost versus spectral radius for the three controllers, where $\rho(A) \in \{0.85, 0.9, 0.95, 1.0, 1.05, 1.1\}$. Curves are averaged over 10 trials, and shaded regions denote 95\% confidence intervals across trials. Left: random bilinear system. Right: multi-block double integrator.}
    \label{fig:rho-sweep}
\end{figure}

Figure~\ref{fig:rho-sweep} shows that the performance gap between \texttt{B-MPC} and the baseline controllers narrows as the dynamics become more stable. When the spectral radius is small, the system is easier to regulate, and the benefit of dual control is reduced. 
The marginal value of improving observability is smaller, so explicitly planning over the belief covariance provides less benefit. Consequently, \texttt{B-MPC} performs similarly to the \texttt{Sep} and \texttt{Sep-MPC}.

\section{Conclusion} \label{sec:conclusion}
We studied quadratic control of linear dynamical systems from bilinear observations, a setting in which the control input affects both state evolution and observation quality, and where the separation principle no longer applies. 
To address this, we proposed a belief-space MPC controller that plans over both the estimated state and the estimation error covariance. 
Our experiments show that \texttt{B-MPC} can outperform both the separation principle controller and an MPC variant of separation control. 
The empirical results further show that the gain of \texttt{B-MPC} is accompanied by low estimation error covariance and by action choices that differ most from the baseline when uncertainty is high, supporting the dual-control interpretation of the method. Overall, these findings suggest that explicitly modeling uncertainty evolution is an important design principle for control under bilinear observations. 
Future work includes establishing rigorous closed-loop stability guarantees for the receding-horizon approximation, developing specialized solvers to efficiently handle the nonconvex planning problem, and extending the policy space to neural-network architectures trained via deep reinforcement learning.

\section{Acknowledgments}
S.D. was partly supported by NSF CCF 2312774, NSF OAC-2311521, NSF IIS-2442137,  a gift to the LinkedIn-Cornell Bowers CIS Strategic Partnership, and an AI2050 Early Career Fellowship program at Schmidt Sciences.

\bibliographystyle{alpha}
\bibliography{arxiv/Bibfiles}

\newcommand{\etalchar}[1]{$^{#1}$}
\begin{thebibliography}{ZRTAm22}

\bibitem[CSJ{\etalchar{+}}25]{choi2025explore}
Sunmook Choi, Yahya Sattar, Yassir Jedra, Maryam Fazel, and Sarah Dean.
\newblock Explore-then-commit for nonstationary linear bandits with latent dynamics.
\newblock {\em arXiv preprint arXiv:2510.16208}, 2025.

\bibitem[GP16]{grune2016nonlinear}
Lars Gr{\"u}ne and J{\"u}rgen Pannek.
\newblock Nonlinear model predictive control.
\newblock In {\em Nonlinear model predictive control: Theory and algorithms}, pages 45--69. Springer, 2016.

\bibitem[Kal60]{kalman1960new}
Rudolph~Emil Kalman.
\newblock A new approach to linear filtering and prediction problems.
\newblock 1960.

\bibitem[KBL{\etalchar{+}}23]{kaufmann2023champion}
Elia Kaufmann, Leonard Bauersfeld, Antonio Loquercio, Matthias M{\"u}ller, Vladlen Koltun, and Davide Scaramuzza.
\newblock Champion-level drone racing using deep reinforcement learning.
\newblock {\em Nature}, 620(7976):982--987, 2023.

\bibitem[KKR22]{kafash2022model}
Sahand~Hadizadeh Kafash, Justin Koeln, and Justin Ruths.
\newblock Model predictive control of bilinear systems as uncertain linear systems.
\newblock In {\em 2022 IEEE Conference on Control Technology and Applications (CCTA)}, pages 562--567. IEEE, 2022.

\bibitem[KLC98]{kaelbling1998planning}
Leslie~Pack Kaelbling, Michael~L. Littman, and Anthony~R. Cassandra.
\newblock Planning and acting in partially observable stochastic domains.
\newblock {\em Artificial Intelligence}, 101(1--2):99--134, 1998.

\bibitem[LK24]{liu2024system}
Diyou Liu and Mohammad Khosravi.
\newblock System identification for linear dynamics with bilinear observation models: An expectation–maximization approach.
\newblock In {\em 2024 IEEE 63rd Conference on Decision and Control (CDC)}, pages 7190--7195, 2024.

\bibitem[May14]{mayne2014model}
David~Q Mayne.
\newblock Model predictive control: Recent developments and future promise.
\newblock {\em Automatica}, 50(12):2967--2986, 2014.

\bibitem[NS24]{nagami2024state}
Keiko Nagami and Mac Schwager.
\newblock State estimation and belief space planning under epistemic uncertainty for learning-based perception systems.
\newblock {\em IEEE Robotics and Automation Letters}, 9(6):5118--5125, 2024.

\bibitem[PT87]{papadimitriou1987complexity}
Christos~H. Papadimitriou and John~N. Tsitsiklis.
\newblock The complexity of markov decision processes.
\newblock {\em Mathematics of Operations Research}, 12(3):441--450, 1987.

\bibitem[PY10]{pardalos2010optimization}
Panos~M Pardalos and Vitaliy~A Yatsenko.
\newblock {\em Optimization and control of bilinear systems: theory, algorithms, and applications}, volume~11.
\newblock Springer Science \& Business Media, 2010.

\bibitem[RMD24]{rawlings2024model}
J.B. Rawlings, D.Q. Mayne, and M.~Diehl.
\newblock {\em Model Predictive Control: Theory, Computation, and Design}.
\newblock Nob Hill Publishing, LLC, 2024.

\bibitem[SCJ{\etalchar{+}}25]{sattar2025sub}
Yahya Sattar, Sunmook Choi, Yassir Jedra, Maryam Fazel, and Sarah Dean.
\newblock Sub-optimality of the separation principle for quadratic control from bilinear observations.
\newblock In {\em 2025 IEEE 64th Conference on Decision and Control (CDC)}, pages 3862--3867. IEEE, 2025.

\bibitem[SJD25]{sattar2025learning}
Yahya Sattar, Yassir Jedra, and Sarah Dean.
\newblock Learning linear dynamics from bilinear observations.
\newblock In {\em 2025 American Control Conference (ACC)}, pages 3109--3115. IEEE, 2025.

\bibitem[XBSA25]{xie2025bilinear}
Yifan Xie, Julian Berberich, Robin Str{\"a}sser, and Frank Allg{\"o}wer.
\newblock Bilinear data-driven min-max mpc: Designing rational controllers via sum-of-squares optimization.
\newblock In {\em 2025 IEEE 64th Conference on Decision and Control (CDC)}, pages 1042--1047. IEEE, 2025.

\bibitem[ZPT23]{zinage2023optimal}
Vrushabh Zinage, Ali~Reza Pedram, and Takashi Tanaka.
\newblock Optimal sampling-based motion planning in gaussian belief space for minimum sensing navigation, 2023.

\bibitem[ZRTAm22]{zheng2022belief}
Dongliang Zheng, Jack Ridderhof, Panagiotis Tsiotras, and Ali-akbar Agha-mohammadi.
\newblock Belief space planning: A covariance steering approach.
\newblock In {\em 2022 International Conference on Robotics and Automation (ICRA)}, pages 11051--11057. IEEE, 2022.

\end{thebibliography}

\newpage
\appendix
\section{Additional Computational Experiments}
\label{app:computational_experiments}

In this appendix, we provide additional experiments from an optimization and implementation perspective. 
The goal is to clarify the computational cost of the proposed belief-space MPC controller and to examine the effect of initialization in the L-BFGS solver used for \texttt{B-MPC}.

\subsection{Wall-clock Runtime} \label{app:wall_clock_runtime}

We compare the wall-clock time of \texttt{Sep}, \texttt{Sep-MPC}, and \texttt{B-MPC} across planning horizons $H$. 
The reported time corresponds to the computation time required to generate one full closed-loop trajectory of length $T=300$, averaged over 10 trials. The system parameter configuration follows Table~\ref{tab:sys-param-config} in Section~\ref{sec:results-best}.

For \texttt{Sep-MPC}, we consider two implementations.
The first implementation solves the finite-horizon deterministic linear-quadratic planning problem through the Riccati recursion, which is described in Section~\ref{sec:controllers}.
This implementation is numerically stable and expected to have nearly the same computational cost as \texttt{Sep}.
In the main experiments, we use the Riccati-based implementation of \texttt{Sep-MPC}.
The second implementation solves the same deterministic planning problem by directly optimizing over the $H$-step planning input sequence using L-BFGS.
This is equivalent to ignoring the covariance component of the belief-space objective in the \texttt{B-MPC} optimization problem while keeping the deterministic state component.

\begin{figure}[ht]
    \centering
    \begin{subfigure}{0.35\textwidth}
        \centering
        \includegraphics[width=\linewidth]{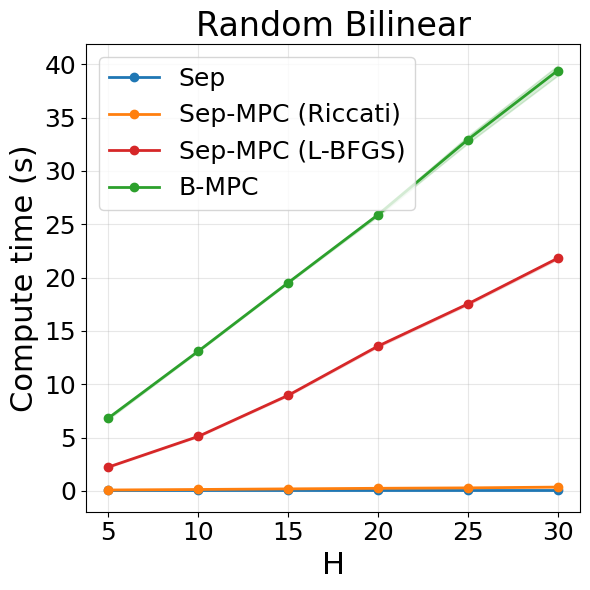}
        \label{fig:h-sweep-time-random}
    \end{subfigure}
    \hspace{10mm}
    \begin{subfigure}{0.35\textwidth}
        \centering
        \includegraphics[width=\linewidth]{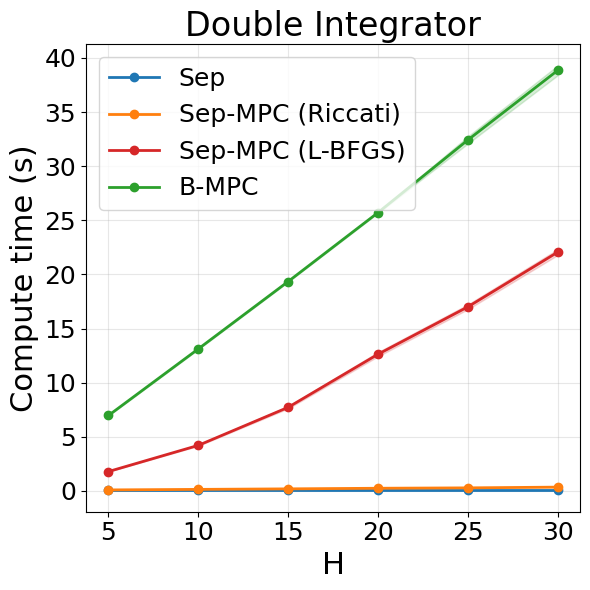}
        \label{fig:h-sweep-time-double}
    \end{subfigure}
    \caption{Computation time versus planning horizon $H$ for the random bilinear and double-integrator systems. The $y$-axis reports wall-clock time in seconds for one full closed-loop trajectory ($T=300$), averaged over 10 trials. Controllers compared are \texttt{Sep}, \texttt{Sep-MPC}, \texttt{B-MPC}. For \texttt{Sep-MPC} controller, it is implemented in two different ways: one is solving the Riccati equation and the other is obtained via L-BFGS optimizer by ignoring the covariance in the belief state.}
    \label{fig:wall_clock_runtime}
\end{figure}

Figure~\ref{fig:wall_clock_runtime} shows the wall-clock comparison.
Although the Riccati-based and L-BFGS-based implementations of \texttt{Sep-MPC} solve the same deterministic planning problem, their runtime difference shows the computational cost of using L-BFGS in place of the finite-horizon Riccati recursion.
The runtime gap between L-BFGS-based \texttt{Sep-MPC} and \texttt{B-MPC} is also substantial, even though the two methods use the same optimization procedure.
The only difference is that \texttt{B-MPC} additionally propagates the nonlinear input-dependent covariance recursion within each objective evaluation.
Thus, the additional computational cost of \texttt{B-MPC} comes from explicitly optimizing over the belief dynamics, especially the covariance evolution.

\subsection{Effect of L-BFGS Initialization} \label{app:lbfgs_initialization}

We next study the effect of initialization on the L-BFGS optimization used by \texttt{B-MPC}.
At each planning step, \texttt{B-MPC} solves a nonconvex finite-horizon optimization problem over the planned control sequence.
A natural question is whether initializing this optimization problem at the \texttt{Sep-MPC} control sequence improves the resulting solution compared to random initialization.
The system parameter configuration is follows Table~\ref{tab:sys-param-config} in Section~\ref{sec:results-best}.

We compare two initialization schemes.
The first is the random initialization used in the main experiments, where the entries of the planned control sequence are initialized independently from $\Ncal(0,1/H)$.
The second is \texttt{Sep-MPC} initialization, where the planned control sequence is initialized using the open-loop solution of \texttt{Sep-MPC} computed from the same estimated state $\hat{x}_{t|t-1}$, without using the covariance component $\Sigma_{t|t-1}$.
We then vary the maximum number of L-BFGS iterations and compare the resulting closed-loop rollout cost.

\begin{figure}[ht]
    \centering
    \begin{subfigure}{0.35\textwidth}
        \centering
        \includegraphics[width=\linewidth]{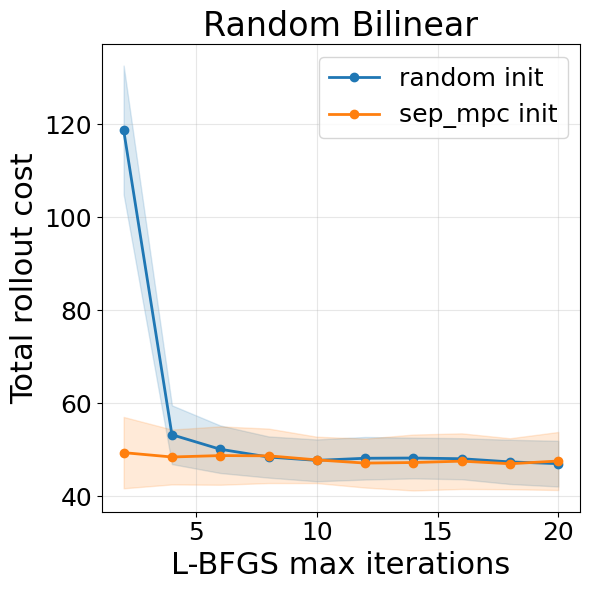}
        \label{fig:l-bfgs-iter-sweep-random}
    \end{subfigure}
    \hspace{10mm}
    \begin{subfigure}{0.35\textwidth}
        \centering
        \includegraphics[width=\linewidth]{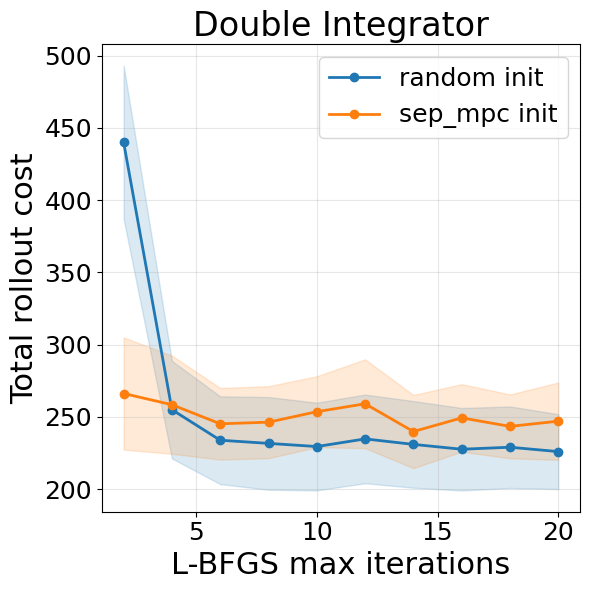}
        \label{fig:l-bfgs-iter-sweep-di}
    \end{subfigure}
    \caption{Total rollout cost versus the maximum number of L-BFGS iterations for \texttt{B-MPC} under two initialization schemes.
    We compare random initialization with \texttt{Sep-MPC} initialization, where the planned control sequence is initialized using the \texttt{Sep-MPC} solution from the same belief state.
    The experiments use planning horizon $H=15$ and trajectory length $T=100$, averaged over 10 trials. The shaded regions denote 95\% confidence intervals across the trials.}
    \label{fig:lbfgs_initialization}
\end{figure}

Figure~\ref{fig:lbfgs_initialization} shows that \texttt{Sep-MPC} initialization does not consistently improve the \texttt{B-MPC} solution.
This effect is especially noticeable in double integrator system, where \texttt{Sep-MPC} initialization can lead to substantially higher cost than random initialization, even after several L-BFGS iterations.
Although \texttt{Sep-MPC} initialization provides a natural initial point, the resulting L-BFGS solution can remain in a more suboptimal region of the nonconvex belief-space objective.
Therefore, we use random initialization in the main experiments.

\end{document}